# STUDY ON OPTIMAL TIMING OF MARK-TO-MARKET FOR CONTINGENT CREDIT RISK CONTROL

Jiali Liao, Drexel University
Theodore V. Theodosopoulos, Drexel University,

**Abstract**

Over-the-counter derivatives have contributed significantly to the effectiveness and efficiency of the international financial system but also entail significant counterparty credit risk. Collateralization is one of the most important and widespread credit risk mitigation techniques used in derivatives transactions. However, the relevant decisions are often made in an ad-hoc manner, without reference to an analytical framework. Very little academic research has addressed the quantitative analysis of collateralization for contingent credit risk control. The issue of mark-to-market timing becomes important for reducing credit exposure of illiquid and long term derivatives contracts due to the difficulty and cost of marking to market. The goal of this research is to propose a framework for minimizing the potential credit exposure of collateralized derivatives transactions by optimizing mark-to-market timing.

## 1. Introduction

OTC derivatives markets have grown exponentially since its introduction in 1980s. The notional amounts of interest rate and currency swaps and interest rate options increased from $865 billion in 1987 to $183 trillion in 2004 and the outstanding notional amount of all OTC derivatives totaled $220 trillion in 2004 (BIS,2004). The rapid growth of OTC derivatives has contributed significantly to the effectiveness and efficiency of the international financial system. OTC derivatives are efficient way to transfer risk and promote optimal risk allocation among agents with different risk tolerance. However, market participants of derivatives transactions faces significant counterparty credit risks, the counterparty may fail to perform on its contingent obligations that are variable and market driven.

The essential underpinning of contingent credit risk management is a meaningful measurement of credit exposure, which is mainly composed of potential credit exposure, a credit risk concept largely unique to OTC derivatives transactions. Potential credit exposure is the amount that one would lose if counterparty default at some future date. The accurate assessment of potential credit exposure is difficult and of great interest to market participants and regulators. This issue becomes more complicated with the rapidly increasing use of collateralization that has proved to be able to significantly reduce credit risk. Collateralization is the most important and widely used methods in practice to mitigate counterparty credit risk. In a typical collateral arrangement, the secured obligation is periodically marked-to-market and the collateral is adjusted to reflect changes in value. The securing party posts additional collateral when the market value has risen, or removes collateral when it has fallen. How to effective measure potential credit exposure of the collateralized portfolio is of significant importance for a wide range of credit risk management applications in derivatives transaction such as capital adequacy, credit limit setting and monitoring etc.

There are three components in credit risk measurement of a financial asset: the probability of default (PD), loss given default (LGD, defined as one minus recovery rate RR), and exposure at default (EAD). While academic literature has focused on the issue of default probabilities and recover rate for pricing credit risk, EAD has been taken for given and lightly touched by the academic models. EAD of loan transactions is typically treated as the book value of assets plus contributions from undrawn credit lines that are calibrated as expectation values by regulatory and propriety models. However, the credit exposure of OTC derivatives transaction is volatile and market driven and the accurate measurement of credit exposure (especially potential credit exposure) is the center feature of effective counterparty risk management. On the other hand, although collateralization has been widely used in practice for reducing credit risk, its impact on credit risk is lightly touched and mainly indirectly studied via its effect on recovery rate and little research exists on how the collateralization affects credit risk exposure and how to determine the major parameters of collateral agreement (such as haircut and margin calls) notwithstanding their importance in practice. Margrabe (1978) mentions the similarity between an exchange option and a margin account and provides the pricing for a very simple marking-to-market. Stulz and Johnson (1985) study the impact of collateralization on the pricing of secured debt using contingent claim analysis. Jokivuolle and Peura (2000) present a model of collateral haircut determination for bank loans. Their model is geared to providing adequate long-to-value ratios, which is similar to the concept of haircut, using a structural credit risk approach. Jokivuolle and Peura (2003) present a model of risky debt in which collateral value is correlated with the possibility of default. The model is then used to study the expected loss given default, primarily as a function of collateral. Aparicio, Felipe M. and Didier Cossin (2003) develop stylized programs providing parties exposed to credit risk with optimal timing and optimal size of the controls required on the collateral they secure, in cases of both complete and partial observation. Cossin and Hricko (2003) present a methodology for haircut determination also using a structural approach but with the final objective of pricing a credit risk instrument backed with collateral. Johannes and Sundaresan (2003) analyze the role of collateral in determining market swap rates. They proved a theory of swap valuation under collateralization and found evidence supporting the presence of costly collateral. They also estimate a terms structure model to characterize the cost of collateral and quantify its effect on swap rates. Cossin and Huang (2003) build a general framework for risk control determination of collateral used in repurchase agreements or repos. They derive how a collateral taker should build collateral policies (i.e. haircut schedules) for the collateral policy to be consistent.

Although marking-to-market secured obligations more frequently has the benefit of reducing the potential exposure, it also incurs substantial operational expenses and the information cost of releasing positions to other dealers, especially for illiquid contract with long maturity. For illiquid OTC derivatives transactions, it is difficult and/or expensive to mark to market the derivatives contracts and collateral held. When the frequency of mark-to-market is constrained, the timing of MTM can provide another way to reduce potential exposure and become important for mitigating credit risk. While the issue of the frequency of mark-to-market and collateral level (Cossin, 2002, 2003) has been studied, no academic research has addressed the quantitative analysis of MTM timing in counterparty potential credit exposure. The goal of this research is to fill in this theoretical gap by providing a general framework that would help evaluate the effect of MTM timing on potential credit exposure of a portfolio. The results of this study have strong implications for academics and practitioners involved in the counterparty risk management of OTC derivatives. The metric of potential credit exposure used in this paper is potential future exposure (PFE), the quantile of maximum credit exposure of all the potential price paths over a specified time interval (such as the lifetime of the contract) at a certain confidence level (95%, 99% etc.). PFE can provide information of worst case over the entire life of a transaction or portfolio. It is a probabilistic metric of future exposure and provides a more accurate and robust analytical metric compared to other metrics like peak exposure or expected exposure.

Dealers in derivatives transactions have a portfolio of counterparties (the firm's portfolio), with each of whom the dealer have a portfolio of derivatives contracts (the counterparty's portfolio). The analysis of potential credit exposure is conducted at three levels of individual transaction, of counterparty and of the firm's portfolio. This paper treats the general case of a single collateralized contract. The simplified case of a single MTM policy is studied and then twice MTM policy are examined. Future research will examine the MTM timing issue for the potential credit exposure at the counterparty level and the firm's portfolio.

This paper is organized as follows. In section 2, we develop the model for single MTM policy and present the numerical results. Section 3 builds the model for twice MTM and displays the numerical results. Section 4 concludes.

**2 Single Mark-to-market Model**

**2.1 Set-up**

This section treats the general case of an individual transaction with unilateral collateral agreement and single marking to market opportunity during the lifetime of the collateralized contract. The decision of optimal MTM timing is modeled as a non-linear optimization problem based on the following assumptions:
(1) The contract has a maturity of $T$ (1, 2, 5 years etc.); the market value of the contract, $V_t$ and $t \in [0,T]$, follows driftless Brownian motion, $BM(0,\sigma)$; the initial market value for the collateral taker is nonnegative, $V_0 \geq 0$.
(2) At contract initiation, the counterparty is required to post a certain amount of collateral ($C_0$, initial margin), which is calculated as a pre-specified percentage ($\beta$) of initial contract value, $C_0 = \beta V_0$. The time to maturity of the collateral is greater than the expiration date of the underlying contract. The market value of collateral ($C_t$) is constant (non-stochastic collateral, i.e. cash).

(3) There is only one marking to market during the lifetime of the contract; there are constant time intervals between candidate MTM days; that is, $\tau = 1,2,...,T-1$. On the pre-specified MTM day ($\tau$), if the MTM value of the contract exceeds the trigger level, a specified percentage of initial market value of contract ($\alpha C_0$), a margin call happens and the counterparty is required to deposit more collateral such that the amount of collateral is brought up to a new level as a percentage of MTM value, that is

If $V_\tau > \alpha C_0$, then $C_\tau = \beta V_\tau$ where $\tau \in (0,T)$

If the MTM value of the contract stays below trigger level, no margin call happens and the collateral remain the same as initial margin, that is,

If $V_\tau \leq \alpha C_0$, then $C_\tau = C_0 = \beta V_0$ where $\tau \in (0,T)$

(4) The maximum credit exposure of the collateralized contract during its lifetime under a specified MTM policy ($E_\tau$) is the maximum amount of the market value of the contract above collateral held during the contract's lifetime. Because the collateral after MTM day may be different than the initial margin, $E_\tau$ is the greater of the maximum credit exposures before MTM day ($E_1$) and that after MTM day ($E_2$). That is,

$$E_\tau = \max\{E_1, E_2\} = \max\left\{\max_{s \in [0,\tau]} V_s - C_0, \max_{s \in [\tau,T]} V_s - C_\tau\right\} = \max\{V_A - C_0, V_B - C_\tau\}$$

where $V_A = \max_{s \in [0,\tau]} V_s$ is the maximum market value of the contract from the contract initiation until MTM day ($[0,\tau]$), which is unobservable but the distribution can be estimated;

$V_B = \max_{s \in [\tau,T]} V_s$ is the maximum market value of the contract from MTM day until maturity ($[\tau,T]$), which is unobservable but the distribution can be estimated;

$V_A$ and $V_B$ are independent conditional on $V_\tau$ according to the Markov property of Brownian motion.

(5) Potential Future Exposure (PFE) is the metric for measuring the potential credit exposure of derivatives transactions. PFE of a collateralized contract is the maximum credit exposure of the collateralized contract expected to occur at a specific confidence level (95%, 99% etc.) over a certain time period (the lifetime of the contract).

Based on the above assumptions, the optimal MTM timing decision is modeled as a nonlinear optimization problem with explicit consideration of the stochastic path properties as below:

$$\underset{\tau}{Min} \quad y$$

$$s.t. \quad f(\tau, y) = q$$

$$y > 0$$

$$\tau \in \{1, 2, ..., T-1\}$$

Where $q = f(\tau, y) = \Pr(E_\tau > y)$ is the probability that $E_\tau$, the maximum credit exposure of a collateralized contract during its lifetime for a given MTM timing policy exceeds a specified risk appetite (a certain level of PFE). The objective function is to minimize PFE at a specified level of confidence level $(1-q)$. Decision variable is MTM time ($\tau$) that can take discrete value among available MTM days between which are constant time intervals.

The function of the probability ($q$) with regard to MTM time ($\tau$) and PFE ($y$) is constructed in the following section. Numerical methods are then used to find the optimal solution of MTM timing ($\tau$) and the corresponding optimal PFE.

**2.2 Model**

This section develops the function of the probability $q$ with regard to MTM time $\tau$, risk appetite $y$ (a specific level of PFE given the volatility rate of the contract) and the parameters of the collateral agreement ($\alpha, \beta$). The market value of the contract is observable only at the initiation, maturity and MTM day ($t = 0, \tau, T$). The collateral after the MTM day is decided by the occurrence of a margin call on that day with regard to the MTM value of the contract ($V_\tau$). Two market scenarios depends on $V_\tau$: (1) Scenario 1 that no margin call happens on MTM day if $V_\tau \leq \alpha C_0$; and (2) Scenario 2 that a margin call happens on MTM day if $V_\tau > \alpha C_0$. Correspondingly, the probability ($q$) is computed in two parts with regard to the occurrence of a margin call on the future MTM day.

$$q = \Pr(E_t > y, V_t \leq \alpha C_0) + \Pr(E_t > y, V_t > \alpha C_0)$$
$$= 1 - \frac{\Pr(V_A - C_0 \leq y, V_\tau \leq \alpha C_0) \cdot \Pr(V_B - C_0 \leq y, V_\tau \leq \alpha C_0)}{\Pr(V_\tau \leq \alpha C_0)} \quad (2.1)$$
$$- \frac{\Pr(V_A - C_0 \leq y, V_\tau > \alpha C_0) \cdot \Pr(V_B - \beta V_t \leq y, V_\tau > \alpha C_0)}{\Pr(V_\tau > \alpha C_0)}$$

Where

$$\Pr(V_\tau \leq \alpha C_0) = \frac{1}{\sqrt{2\pi\tau}\sigma} \int_{-\infty}^{\alpha C_0} \exp\left(-\frac{(x-V_0)^2}{\sigma^2 \tau}\right) dx = 1 - \frac{1}{2} erfc\left(\frac{\alpha C_0 - V_0}{\sigma\sqrt{2\tau}}\right) \quad (2.2)$$

Equation 2.2 calculates the probability of margin call not triggered the MTM value of the contract on MTM day is below the trigger level, that is $V_\tau \leq \alpha C_0$. The MTM value of the contract ($V_\tau$) follows a normal distribution with mean of 0 and variance of $\sigma^2 \tau$, that is $V_\tau \sim N(0, \sigma^2 \tau)$.

$$\Pr(V_A - C_0 \leq y, V_\tau \leq \alpha C_0)$$
$$= \Pr(V_\tau \leq \alpha C_0) - \Pr(V_A - C_0 > y | V_t \leq \alpha C_0) \cdot \Pr(V_t \leq \alpha C_0) \quad (2.3)$$
$$= 1 - \frac{1}{2} erfc\left(\frac{\alpha C_0 - V_0}{\sigma\sqrt{2\tau}}\right) - \left[1 - \frac{1}{2} erfc(\frac{\alpha C_0 - (2y + 2C_0 - V_0)}{\sigma\sqrt{2\tau}})\right]$$

Equation 2.3 computes the probability of the maximum credit exposure before MTM timing below the risk appetite under scenario 1. The distribution of $V_A$, the maximum market value of the contract from the contract initiation until MTM day, conditional on the MTM value of the contract on MTM day $(V_\tau = x \leq \alpha C_0)$ can be estimated using the distribution of the maximum of a Brownian Bridge with given starting and ending values $(V_0, x)$. The proof sees Appendix A.

$$\Pr(V_B - C_0 \leq y, V_\tau \leq \alpha C_0)$$
$$= 1 - \frac{1}{2} erfc\left(\frac{\alpha C_0 - V_0}{\sigma\sqrt{2\tau}}\right) - \int_{-\infty}^{\alpha C_0}\left[erfc\left(\frac{C_0 + y - x}{\sqrt{2(T-\tau)}\sigma}\right) * \frac{1}{\sqrt{2\pi\tau}\sigma} e^{-\frac{(x-V_0)^2}{2\sigma^2\tau}}\right]dx \quad (2.4)$$

Equation 2.4 computes the probability that maximum credit exposure after MTM day until maturity stays below risk appetite with no margin call on MTM day. The market value of the contract during $[\tau, T]$ follows Brownian motion $BM(0, \sigma^2)$ starting from $V_\tau$. Using reflection principle of Brownian motion, we can easily derive the conditional probability of $V_B$ given $V_\tau$.

$$\Pr(V_A - C_0 \leq y, V_\tau > \alpha C_0)$$
$$= \frac{1}{2} erfc\left(\frac{\alpha C_0 - V_0}{\sigma\sqrt{2\tau}}\right) - \frac{1}{2} erfc\left(\frac{\alpha C_0 - (2y + 2C_0 - V_0)}{\sigma\sqrt{2\tau}}\right) + 1 - erfc\left(\frac{y + C_0 - V_0}{\sigma\sqrt{2\tau}}\right) \quad (2.5)$$

Equation 2.5 computes the probability of maximum credit exposure from the contract initiation to MTM day below the risk appetite with margin call on MTM day. The derivation of the conditional distribution of $V_A$ conditional on $V_\tau$ characterized by Brownian Bridge sees the Appendix.

$$\Pr(V_B - \beta V_t \leq y, V_\tau > \alpha C_0)$$
$$= \frac{1}{2} erfc\left(\frac{\alpha C_0 - V_0}{\sigma\sqrt{2\tau}}\right) - \int_{\alpha C}^{-\infty}\left[erfc\left(\frac{\beta x + y - x}{\sqrt{2(T-\tau)}\sigma}\right) * \frac{1}{\sqrt{2\pi\tau}\sigma} e^{-\frac{(x-V_0)^2}{2\sigma^2\tau}}\right]dx \quad (2.6)$$

Equation 2.6 computes the probability of maximum credit exposure after MTM day until maturity below risk appetite with margin call on MTM day.

**2.3 Numerical Examples**

After the construction of the function of the probability with regard to risk appetite and MTM time and the parameters, this section gives numerical examples and uses numerical methods to solve the optimal MTM time for minimizing PFE.

The benchmark case uses the following values of the main parameters in this model, $q = 0.05, \alpha = 0.9, \beta = 1.1, V_0 = 1, T = 24, \sigma = 0.2$. This case studies the optimal MTM time to minimize PFE at 95% confidence level, starting with the market value of the underlying contract ($V_0$) at initiation is $1 million, when the collateral pledged is 110% of $V_0$. On MTM day, if market value of the underlying contract exceeds 90% of collateral pledged, margin call is triggered and the counterparty is asked to pose additional collateral such that the new collateral is 110% of the observed MTM value of the underlying contract; otherwise, no margin call is triggered and no action is needed. The market value of the underlying contract follows driftless Brownian motion with monthly volatility 20%. Figure 2.1 displays that there exist optimal time of marking to market for minimizing PFE and the approximate optimal solution is around $10^{th}$ month for the contract with a maturity of 24 months.

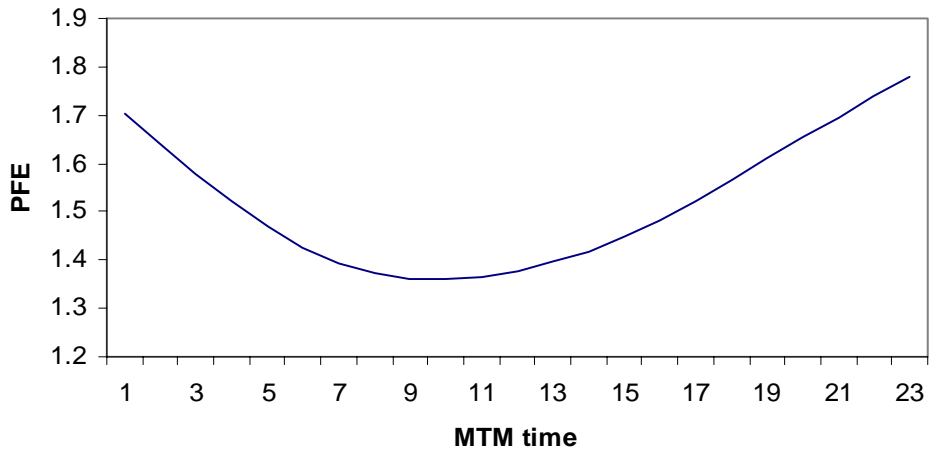

Figure 2.1 PFE w.r.t MTM time for the benchmark case

Figure 2.2 presents the pattern of the probability $q$ with regard to PFE and MTM day ($\tau$) for contract with different volatilities using the analytic formula for a contract of 24 months maturity. Figure 2.3 displays the results using Monte Carlo Simulation methods.

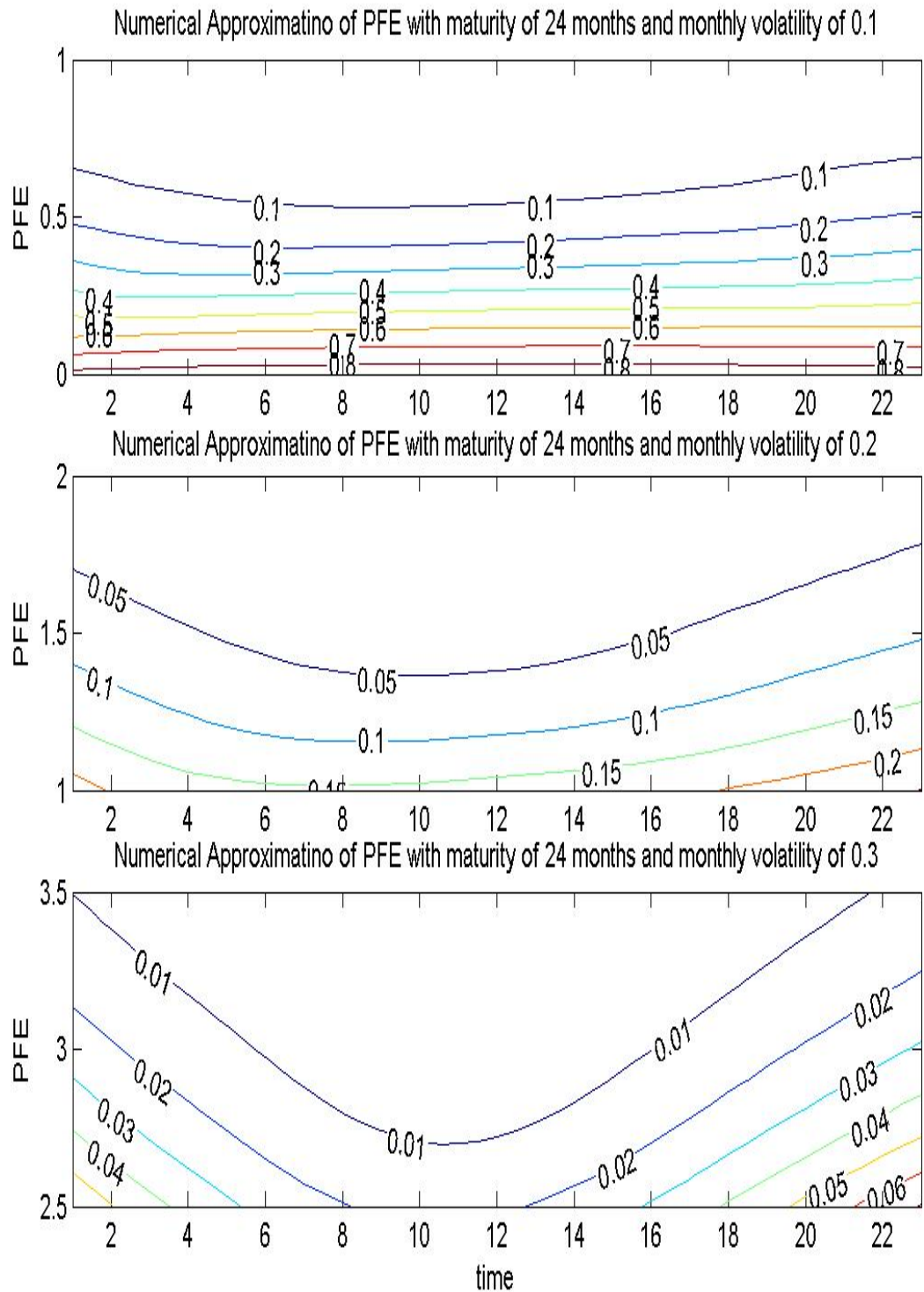

Figure 2.2 Numerical examples with maturity of 24 months

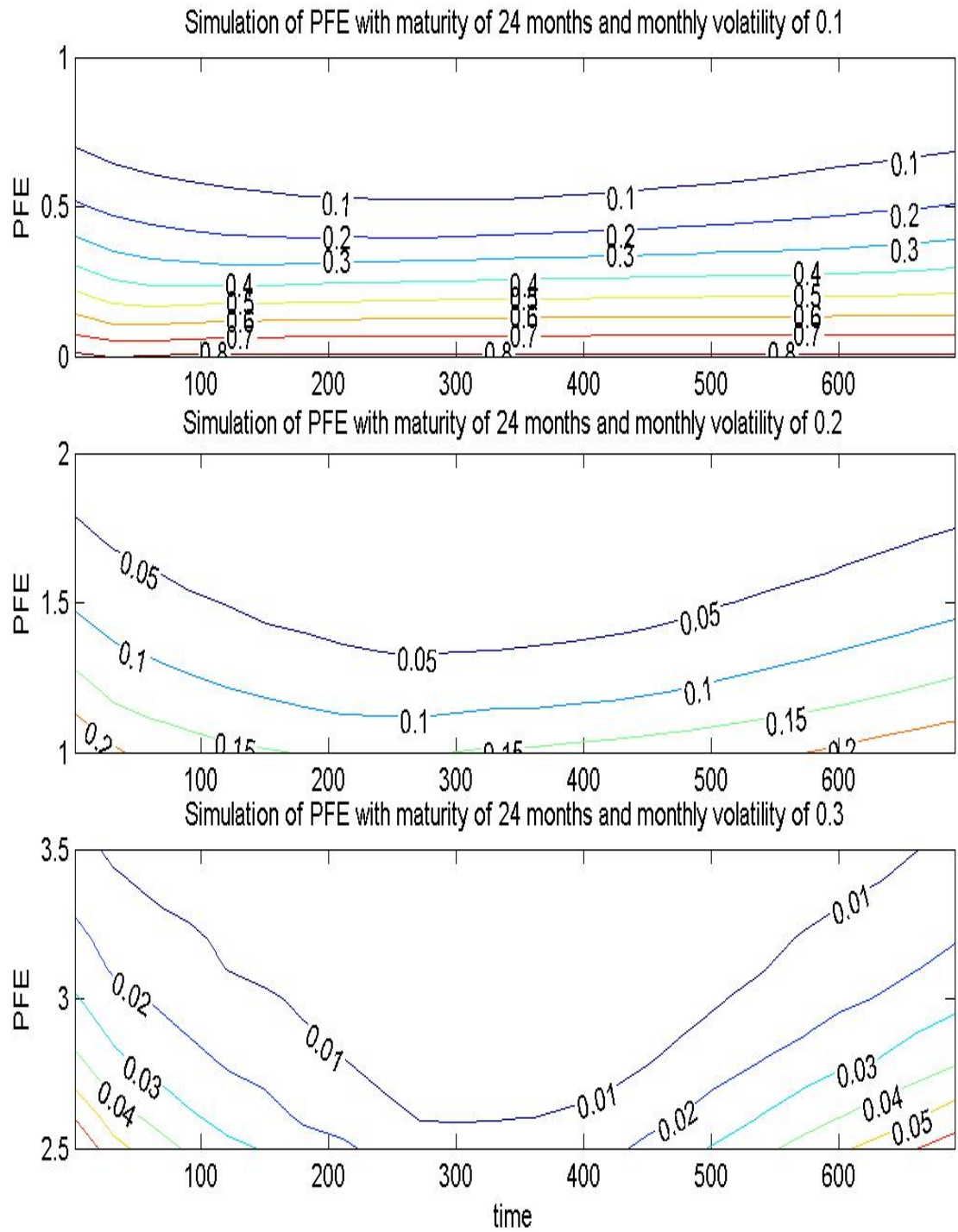

Figure 2.3 Simulation of q w.r.t PFE and MTM time with maturity of 24 months

With the base case as a benchmark, the effect on optimal MTM timing and minimum PFE of changing the values of the different parameters of the models is investigated.

(1) Changes in $T$ (maturity of the underlying contract)

| T | MTM TIME | PFE |
|---|---|---|
| 12 | 5 | 0.9325 |
| **24** | **10** | **1.3602** |
| 36 | 15 | 1.6884 |

As expected, the longer the time to maturity of the underlying contract, the higher the PFE. The optimal MTM time is around the half way (0.4166) of the maturity.

(2) Change in $q$ (the probability of the maximum credit exposure above risk appetite)

| Q | MTM TIME | PFE |
|---|---|---|
| 0.1 | 9 | 1.1515 |
| **0.05** | **10** | **1.3602** |
| 0.01 | 11 | 1.7697 |

Also as expected, the smaller the confidence level, the higher the PFE. The optimal MTM time is closer to the maturity.

(3) Change in $\sigma$ (the monthly volatility of the underlying contract)

| $\sigma$ | MTM TIME | PFE |
|---|---|---|
| 0.1 | 10 | 0.4163 |
| **0.2** | **10** | **1.3602** |
| 0.3 | 10 | 2.0903 |

The higher the volatility of the underlying contract, the higher PFE; however, the optimal MTM time is not affected by the volatility.

(4) Change in $V_0$ (initial contract value)

| V0 | MTM TIME | PFE |
|---|---|---|
| 0 | 10 | 1.4602 |
| 0.5 | 10 | 1.4102 |
| **1** | **10** | **1.3602** |
| 1.5 | 10 | 1.3102 |
| 2 | 10 | 1.2603 |

The higher the initial contract value, the smaller PFE; the optimal MTM time remains the same for different initial contract value.

(5) Change in $\alpha$ (the multiplier for trigger level of margin call)

| $\alpha$ | MTM TIME | PFE |
|---|---|---|
| 0.5 | 11 | 1.3966 |
| 0.8 | 10 | 1.3613 |
| **0.9** | **10** | **1.3602** |

| 1 | 10 | 1.3636 |

There is a value of alpha with which PFE is minimal.

(6) Change in $\beta$ (multiplier for computing required collateral)

| BETA | MTM TIME | PFE |
|------|----------|--------|
| 1    | 10       | 1.4817 |
| 1.1  | 10       | 1.3602 |
| 1.5  | 10       | 0.9703 |
| 1.7  | 11       | 0.8321 |
| 1.9  | 12       | 0.7026 |
| 2    | 13       | 0.6364 |

As expected, the higher the required collateral level, the smaller the PFE. With the increase of $\beta$, the optimal MTM time moves closer the maturity.

## 3 Optimal Timing Analysis of Twice Mark-to-market Policy

This section extends to the twice MTM policy where the timing decision is made simultaneously at the contract initiation or in a sequential manner. Under the simultaneous decision scheme, the collateral taker decides at the contract initiation the two days when to mark the contract to the market; while under the sequential decision scheme, the collateral taker decides the first MTM time at the beginning of the contract and set the second MTM time on the first day when the contract is marked to market

### 3.1 Twice MTM of Simultaneous Timing Decision

This section treats the general case of an individual transaction with unilateral collateral agreement under twice MTM policy with simultaneous timing decision. The assumptions are the same as those for single MTM policy except that there are two marking to market opportunities during the lifetime of the contract.

The collateral taker needs to decide the two MTM days $(\tau_1, \tau_2)$ at the beginning of the contract; there are constant time intervals between candidate MTM days; that is, $\tau_1 \in \{1,2,...,T-1\}$ and $\tau_2 \in \{\tau_1+1,2,...,T-1\}$. On the first pre-specified MTM day $(\tau_1)$, if the market value of the contract exceeds the trigger level which is set by a percentage $(\alpha)$ of collateral posted, a margin call happens and the counterparty will be required to deposit more collateral; otherwise, no more collateral is required. That is,

$$C_1 = \begin{cases} \beta V_0 & V_1 \leq \alpha C_0 \\ \beta V_1 & V_1 > \alpha C_0 \end{cases} \text{ if }$$

where $V_1$ is the market value of the underlying contract at $\tau_1$

On the second pre-decided MTM day $(\tau_2)$, if the market value of the underlying contract exceeds the trigger level which is set by a multiple of collateral posted, a margin call happens and the counterparty will be required to deposit more collateral; otherwise, no more collateral is required. That is

$$C_2 = \begin{cases} \beta V_0 & V_1 \leq \alpha C_0 \text{ and } V_2 \leq \alpha C_0 \\ \beta V_1 & \text{if } V_1 > \alpha C_0 \text{ and } V_2 \leq \alpha C_1 = \alpha \beta V_1 \\ \beta V_2 & V_1 > \alpha C_0 \text{ and } V_2 > \alpha C_1 = \alpha \beta V_1 \end{cases}$$

where $V_2$ is the market value of the underlying contract at $\tau_2$.

Under a given MTM policy with specified MTM days $(\tau_1, \tau_2)$, the maximum credit exposure of a collateralized contract during its lifetime is the maximum amount of the market value of the contract above collateral held. Because the collateral after MTM day may be different than that before the MTM day and there are two MTM days, $E_\tau$ is the greatest of the maximum credit exposures between the contract initiation and the first MTM day ($E_1$), the maximum credit exposure between the first and the second MTM day ($E_2$) and the maximum credit exposure after the second MTM day until maturity ($E_3$).

$E_\tau = \max\{E_1, E_2, E_3\}$, where $E_1, E_2$ and $E_3$ is the maximum credit exposure during time periods $[0, \tau_1]$, $[\tau_1, \tau_2]$ and $[\tau_2, T]$ respectively.

$E_1 = \max_{s \in [0, \tau_1]} V_s - C_0 = V_A - C_0$, where $V_A$ is the maximum market value of the underlying contract during $[0, \tau_1]$

$E_2 = \max_{s \in [\tau_1, \tau_2]} V_s - C_1 = V_B - C_1$, where $V_B$ is the maximum market value of the underlying contract during $[\tau_1, \tau_2]$

$E_3 = \max_{s \in [\tau_2, T]} V_s - C_2 = V_C - C_2$, where $V_C$ is the maximum market value of the underlying contract during $[\tau_2, T]$

The twice MTM time decision is models as nonlinear optimization problem with explicit consideration of the stochastic path properties:

$$\begin{aligned} \min_{\tau_1, \tau_2} \quad & y \\ \text{s.t.} \quad & f(\tau, y) = q \\ & y > 0 \\ & \tau_1 \in \{1, 2, ..., T-2\} \\ & \tau_2 \in \{\tau_1 + 1, 2, ..., T-1\} \end{aligned} \quad (3.1)$$

Where $q = f(\tau, y) = \Pr(E_\tau > y)$ is the probability that maximum credit exposure of the collateralized contract during its lifetime exceeds the risk appetite (a specific value of the PFE).

The objective function is to minimize PFE at a specified level of confidence level $(1-q)$, 95%, 99% etc. The decision variable is MTM time $(\tau_1, \tau_2)$ that can take discrete values from available MTM days between which are constant time intervals.

The framework first derives the function of the probability ($q$) with regard to MTM time ($\tau$) and PFE ($y$) and then uses numerical methods to identify the optimal

MTM timing solution and the corresponding optimal PFE that will be compared to those of sequential timing decision.

The probability is computed in four parts corresponding to four scenarios of the likelihood of margin calls on the two pre-specified MTM days:

(1) No margin call happens on either of the two MTM days if $V_1 \leq \alpha C_0$ and $V_2 \leq \alpha C_1$, where $C_2 = C_1 = C_0 = \beta V_0$.

(2) No margin call happens on first MTM day but a margin call happens on second MTM day if $V_1 \leq \alpha C_0$ and $V_2 > \alpha C_1$, where $C_1 = C_0 = \beta V_0$ and $C_2 = \beta V_1$.

(3) A margin call happens on first MTM day but no margin call happens on second MTM day if $V_1 > \alpha C_0$ and $V_2 \leq \alpha C_1$, where $C_1 = \beta V_1$ and $C_2 = C_1 = \beta V_1$.

(4) A margin call happens on both MTM days if $V_1 > \alpha C_0$ and $V_2 > \alpha C_1$, where $C_1 = \beta V_1$ and $C_2 = \beta V_2$.

## 3.2 Twice Mark-to-market with sequential timing decision

This section treats the general case of an individual transaction with unilateral collateral agreement under twice MTM policy with sequential timing decision. The assumptions are the same as the single MTM model except that there are two MTM opportunities during the lifetime of the contract.

The PFE before the first MTM day is equal to the greater of the quantile (risk appetite) at a confidence level of $(1-q)$ and zero;

$$y_A = \max\{f^{-1}(q_A, x, \tau_1), 0\}$$

and

$$q_A = f(y_A, \tau_1, x) = \Pr(E_A > y_A \mid V_1 = x)$$
$$= \Pr\left(\max_{s \in [0,\tau_1]} V_s - C_0 > y_A \mid V_1 = x\right)$$
$$= \Pr(V_A - C_0 > y_A \mid V_1 = x)$$

$q_A$ is the conditional probability that $E_A$ (the maximum credit exposure between the contract initiation and the first MTM day $[0, \tau_1]$) exceeds a risk appetite given the MTM value of the contract; and $V_A = \max_{s \in [0,\tau_1]} V_s$ is the maximum market value of the contract between the contract initiation and the first MTM day, whose value is unobservable but distribution can be estimated.

The PFE after the first MTM day ($y_B$) is the minimal PFE calculated using single MTM timing decision model with initial contract value $V_1$, the initial margin of $C_1$, the maturity of $T - \tau_1$ and the same margin call policy set before.

$$\min_{\tau_2} \; y_B$$

$$s.t.$$

$$g(y_B, \tau_2, x) = q_B$$

$$y_B \geq 0$$

$$\tau_2 \in \{\tau_1 + 1, \tau_1 + 2, ..., T - 1\}$$

where

$$q_B = g(y_B, x, \tau_2) = \Pr(E_B > y \mid V_1 = x)$$

$$= \Pr\left(\max\left\{\max_{s \in [\tau_1, \tau_2]} V_s - C_1, \max_{s \in [\tau_2, T]} V_s - C_2\right\} > y \mid V_1 = x\right)$$

$$= \Pr(\max\{V_B - C_1, V_C - C_2\} > y \mid V_1 = x)$$

$q_B$ is the conditional probability that $E_B$ (the maximum credit exposure of the collateralized contract between the first MTM day and the maturity day $[\tau_1, T]$) exceeds a risk appetite given the MTM value of the contract on first MTM day $(V_1 = x)$; $E_B$ is the greater of the maximum credit exposure between the first and the second MTM day ($\max_{s \in [\tau_1, \tau_2]} V_s - C_1$) and that between the second MTM day and the maturity day ($\max_{s \in [\tau_2, T]} V_s - C_2$); $V_B$ is the maximum market value of the contract during $[0, \tau_1]$ and $V_C$ is the maximum market value of the contract during $[\tau_2, T]$, both of which are unobservable but the distribution can be estimated.

Based on the above assumptions, the sequential timing decision of twice MTM policy is models as two-stage stochastic dynamic programming problem.

$$\min_{\tau_1, \tau_2(\tau_1)} \; E\, y$$

$$s.t.$$

$$y = \max\{y_A, y_B\}$$

$$\tau_1 \in \{1, 2, ..., T - 2\} \quad (3.2)$$

$$\tau_2 = \arg\min_{\tau_2 \in \{\tau_1 + 1, ..., T - 1\}} y_B$$

The objective function is to minimize the expected PFE over the distribution of the MTM value of the contract on the first MTM day. The PFE of contract during its lifetime conditional on the MTM value of the contract is the greater of the PFE before the MTM value and after the MTM value defined above. The function of the expected PFE with regard to the twice MTM time is developed and then numerical methods are used to find the optimal MTM days and the results are compared with the simultaneous timing decision.

The expected PFE is computed based on the distribution of the MTM value of the contract on the first MTM day. The MTM value of the contract on the first MTM day ($V_1$) follows a normal distribution with a mean of zero and a variance of $\sigma^2 \tau_1$, that is $V_1 \sim N(0, \sigma^2 \tau_1)$. $V_1$ defines two scenarios: (1) Scenarios 1 that no margin call happens if $V_1 \leq \alpha C_0$; (2) Scenario 2 that a margin call may occur if $V_1 > \alpha C_0$. The occurrence of a margin call on the first MTM day decides the collateral level for the second stage timing decision is dependent on; therefore, the expected value of PFE is computed in two parts with regard to $V_1$.

$$Ey = \frac{1}{\sqrt{2\pi\tau_1}\sigma} \int_{-\infty}^{+\infty} y \bullet \exp\left(-\frac{(x-V_0)^2}{2\sigma^2\tau_1}\right) dx$$

$$= \frac{1}{\sqrt{2\pi\tau_1}\sigma} \int_{-\infty}^{+\infty} \max\{y_A, y_B\} \bullet \exp\left(-\frac{(x-V_0)^2}{2\sigma^2\tau_1}\right) dx \qquad (3.3)$$

$$= \frac{1}{\sqrt{2\pi\tau_1}\sigma} \int_{-\infty}^{\alpha C_0} y_1 \bullet \exp\left(-\frac{(x-V_0)^2}{2\sigma^2\tau_1}\right) dx + \frac{1}{\sqrt{2\pi\tau_1}\sigma} \int_{\alpha C_0}^{\infty} y_2 \bullet \exp\left(-\frac{(x-V_0)^2}{2\sigma^2\tau_1}\right) dx$$

$y_1$ is a specific value of the PFE of the collateralized contract during its lifetime under scenario 1 defined by $V_1 = x \leq \alpha C_0$ such that no margin call happens on the first MTM day ($\tau_1$). It is the greater of the PFE between the contract initiation and the first MTM day ($[0,\tau_1]$) and the PFE between the first MTM day and the maturity day ($[\tau_1,T]$). Thus,

$$y_1 = \max\{y_A^1, y_B^1 *\}$$

where $y_A^1$ is the PFE during $[0,\tau_1]$ given $V_1 = x \leq \alpha C_0$ and $y_B^1$ is the optimal PFE during $[\tau_1,T]$ given $V_1 = x \leq \alpha C_0$ solved by the single MTM timing decision model.

$$y_A^1 = f_1^{-1}(q_A^1, \tau_1)$$

where

$$q_A^1 = f_1(y_A^1, \tau_1) = \Pr(E_A > y | V_1 \leq \alpha C_0)$$
$$= \Pr(V_A - C_0 > y | V_1 = x \leq \alpha C_0) \qquad (3.4)$$
$$= \exp\left(-\frac{2(y+C_0)^2 - 2(y+C_0)(x+V_0) + 2xV_0}{\sigma^2\tau_1}\right)$$

$q_A^1$ describes the conditional probability of $E_A$ (the maximum credit exposure of the collateralized contract during $[0,\tau_1]$) exceeding the risk appetite given $V_1 = x \leq \alpha C_0$. The conditional distribution of $E_A$ given $V_1 = x \leq \alpha C_0$ is characterized by the distribution of the maximum of a Brownian Bridge with given beginning and ending value of $(V_0, x)$. Numerical method is then used to find out $y_A^1$ at a specific level of $q_A^1$.

$y_B^1 *$ is the optimal PFE during $[\tau_1,T]$ under the optimal timing decision solved using single MTM timing model given $V_1 = x \leq \alpha C_0$ such that no margin call happens on the first MTM day.

$$\min_{\tau_2} \quad y_B^1$$
$$s.t. \quad g_1(y_B^1, \tau_2, V_1) = q_B^1 \quad (3.5)$$
$$y_B^1 \geq 0$$
$$\tau_2 \in \{\tau_1 + 1, \tau_1 + 2, ..., T - 1\}$$

where
$$q_B^1 = g_1(y_B^1, V_1, \tau_2)$$
$$= \Pr(E_B > y \mid V_1 \leq \alpha C_0)$$
$$= \Pr(\max\{V_B - C_1, V_C - C_2\} > y \mid V_1 \leq \alpha C_0)$$

$y_2$ is a specific value of the PFE of the collateralized contract during its lifetime under scenario 2 defined by $V_1 = x > \alpha C_0$ such that a margin call happens on the first MTM day. It is the greater of the PFE between the contract initiation and the first MTM day ($[0, \tau_1]$) and the PFE between the first MTM day and the maturity day ($[\tau_1, T]$). Thus,
$$y_2 = \max\{y_A^2, y_B^2 *\}$$
where $y_A^2$ is the PFE during $[0, \tau_1]$ given $V_1 = x > \alpha C_0$ and $y_B^2$ is the optimal PFE during $[\tau_1, T]$ given $V_1 = x > \alpha C_0$ solved by the single MTM timing decision model.
$$y_A^2 = f_2^{-1}(q_A^2, \tau_1), \quad \text{where} \quad q_A^2 = f_2(y_A^2, \tau_1) = \Pr(E_A > y \mid V_1 > \alpha C_0)$$

$q_A^2$ describes the conditional probability of $E_A$ (the maximum credit exposure during $[0, \tau_1]$) exceeding a risk appetite given $V_1 = x > \alpha C_0$. The conditional distribution of $E_A$ given $V_1 = x \leq \alpha C_0$ is characterized by the distribution of the maximum of a Brownian Bridge with given beginning and ending value of $(V_0, x)$. Numerical method is then used to find out $y_A^2$ at a specific level of $q_A^2$.

$y_B^2 *$ is the optimal PFE during $[\tau_1, T]$ under the optimal timing decision solved using single MTM timing model given $V_1 = x > \alpha C_0$ such that no margin call happens on the first MTM day.
$$\min_{\tau_2} \quad y_B^2$$
$$s.t. \quad g_2(y_B^2, \tau_2, x) = q_B^2 \quad (3.6)$$
$$y_B^2 \geq 0$$
$$\tau_2 \in \{\tau_1 + 1, \tau_1 + 2, ..., T - 1\}$$
where
$$q_B^2 = g_2(y_B^2, x, \tau_2)$$
$$= \Pr(E_B > y \mid V_1 = x > \alpha C_0)$$
$$= \Pr(\max\{V_B - C_1, V_C - C_2\} > y \mid V_1 = x > \alpha C_0)$$

## 3.3 Numerical Results

After the construction of the function of the probability under simultaneous and sequential timing decision with regard to risk appetite and MTM time and the parameters defining the contract value and collateralization, this section gives numerical examples and uses numerical methods to solve the optimal MTM time for minimizing PFE.

The benchmark case uses the following values of the main parameters in this model: $q = 0.05, \alpha = 0.9, \beta = 1.1, V_0 = 1, T = 12, \sigma = 0.1$. This case studies the optimal time for making to market in order to minimize PFE at 95% confidence level, starting with the market value of the contract at initiation ($V_0$) is $1 million, when the collateral pledged is 110% of $V_0$. On the MTM days, if the market value of the contract exceeds 90% of collateral pledged, margin call is triggered and the counterparty is asked to pose additional collateral such that the new collateral is 110% of the MTM value of the contract; otherwise, no margin call is triggered and no action is needed. The market value of the contract follows driftless Brownian motion with monthly volatility 10%. Table 3.1 compares the results of optimal MTM timing and the corresponding optimal PFE under single MTM policy, twice MTM with simultaneous timing decision and twice MTM with sequential timing decision schemes. It is shown that there exists optimal time of marking to markets for minimizing PFE under all the MTM policies. The optimal PFE under twice MTM with sequential timing decision is minimal among the three MTM policies. The first MTM time under twice MTM policy is the same and earlier than the optimal MTM time under single MTM policy because of the existence of the second MTM opportunity.

Table 3.1 Comparison of PFE with different time decisions

| Single MTM | | Twice MTM (simultaneous) | | | Twice MTM (sequential) | |
|---|---|---|---|---|---|---|
| T | PFE | t1 | t2 | PFE | t1 | PFE |
| 1 | 0.5144 | 1 | 6 | 0.3726 | 1 | 0.3667 |
| 2 | 0.4736 | 2 | 6 | 0.3443 | 2 | 0.3302 |
| 3 | 0.4397 | 3 | 7 | 0.3252 | 3 | 0.2964 |
| 4 | 0.4206 | **4** | **8** | **0.3189** | **4** | **0.2822** |
| **5** | **0.4163** | 5 | 8 | 0.3191 | 5 | 0.2918 |
| 6 | 0.4224 | 6 | 9 | 0.3272 | 6 | 0.3031 |
| 7 | 0.4369 | 7 | 9 | 0.3443 | 7 | 0.3211 |
| 8 | 0.4595 | 8 | 10 | 0.3666 | 8 | 0.3493 |
| 9 | 0.4886 | 9 | 10 | 0.3938 | 9 | 0.3764 |
| 10 | 0.5198 | 10 | 11 | 0.4202 | 10 | 0.402 |
| $q = 0.05, T = 12, \sigma = 0.1$ | | | | | | |

Figure 3.1 to 3.3 presents the PFE with respect to MTM time under three different marking to market policies for different maturities of the contract.

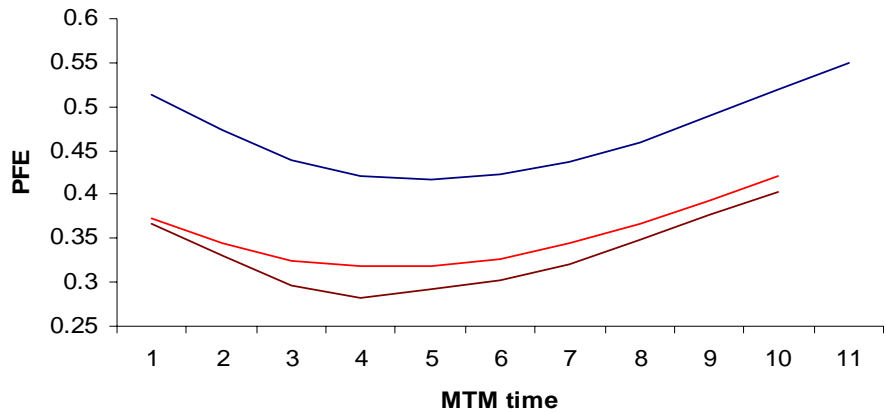

Figure 3.1 PFE w.r.t first MTM time with T=12 and $\sigma = 0.1$

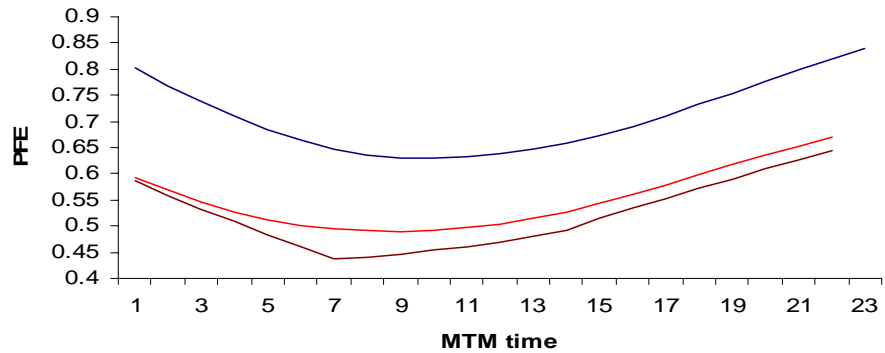

Figure 3.2 PFE w.r.t first MTM time with T=24 and $\sigma = 0.1$

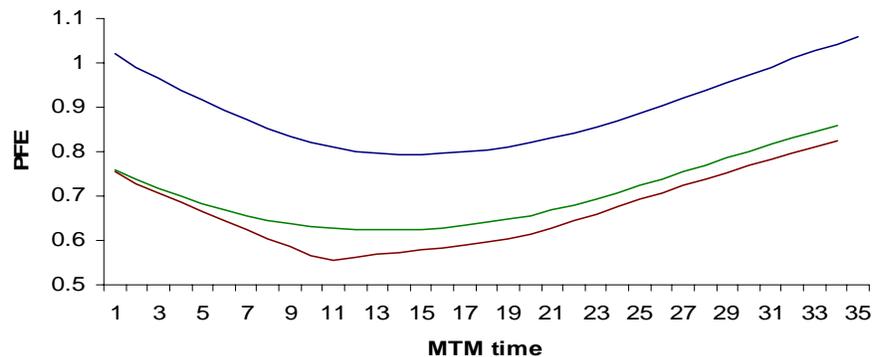

Figure 3.3 PFE w.r.t first MTM time with T=36 and $\sigma = 0.1$

As expected, twice MTM time decisions reduce PFE compared to single MTM time decision. And sequential MTM time decision does better in minimizing PFE than simultaneous decision because the sequential MTM timing utilize more information about market value than simultaneous decision. The optimal first MTM time under sequential time decisions is earlier than that under simultaneous time decisions when maturity increases. The longer the lifetime of the underlying contract, the higher degree the sequential twice MTM time decision reduces PFE.

## 4 Conclusions

In this paper, we have developed a general model of study the optimal marking to market timing in collateralization to minimize the potential future exposure of a single collateralized contract. The simplified case of a single marking to market policy is studied first and then twice marking to market policy is examined. It is proved that there exists an optimal MTM time for minimizing PFE under single and twice marking to market policies with timing decision made in simultaneous and sequential manner. Twice MTM policies are better in reducing PFE than single MTM policy because the additional MTM opportunity. The twice MTM of sequential timing decision is better than that of simultaneous timing decision because the second stage timing decision is an optimal decision based on future market scenario utilizing new upcoming information. This study has important policy implication for contingent credit risk control and lays the foundation for the further analysis of the aggregate contingent credit exposure of derivatives portfolio taking into account the netting effect.

**Appendix: Proof of Distribution of the Maximum of Brownian Bridge**

Given $V_0$ and $V_t = x$, $V_s$ ($s \in [0,t]$) is independent of $V_t$.

For any $\varepsilon > 0$, let $\Pr_{b,\varepsilon} = \Pr(\sup_{s \in [0,t]} V_s > b$ for some $s \in [0,t] \| V_t - x | < \varepsilon)$, $x < b$

Let $\tau = \min\{s : V_s = b, s \in [0,t]\}$, $\tau$ is stopping time. Let $B = \{\tau < t\}$. $V_\tau = b$ on B by sample continuity. By the sample continuity and intermediate value theorem, $\Pr(B) > 0$. For each measurable event A, let $\Pr_b(A) = \Pr(A | B)$. Then with respect to $\Pr_b$, the process $W_h = V_{\tau+h} - V_\tau$ has the law of Brownian motion and is independent of $B^B_{\tau+} = \{D \cap B : D \in B_{\tau+}\}$, by the strong Markov property. Let $h = t - \tau$. Then $h$ is a $B^B_{\tau+}$ measurable random variable, so it is independent of $W_h$ for $\Pr_b$.

$$\Pr_{b,\varepsilon} = \frac{\Pr(\tau < t, |W_h + b - x| < \varepsilon)}{\Pr(|V_t - x| < \varepsilon)} = \frac{\Pr(\tau < t) \cdot \Pr_b(|W_h + b - x| < \varepsilon)}{\Pr(|V_t - x| < \varepsilon)}, \quad 0 < \varepsilon < b$$

The distribution of $W_h$ and $-W_h$ for $\Pr_b$ are equal, thus

$$\Pr_b(|W_h + b - x| < \varepsilon) = \Pr_b(|-W_h - b + x| < \varepsilon) = \Pr_b(|W_h - b + x| < \varepsilon)$$

By definition of $W_s$ and $V_\tau = b$, $W_h - b + x = V_t - 2b + x$ whenever $\tau < t$. By sample continuity and since $0 < \varepsilon < b$, $|V_t - 2b + x| < \varepsilon$ implies $V_t > b$, so $\tau < t$. It follows that

$$\Pr_b(|W_h + b - x| < \varepsilon) = \frac{\Pr(|W_h + b - x| < \varepsilon)}{\Pr(\tau < t)} = \frac{\Pr(|W_h - b + x| < \varepsilon)}{\Pr(\tau < t)} = \frac{\Pr(|V_t - 2b + x| < \varepsilon)}{\Pr(\tau < t)},$$

then

$$\Pr_{b,\varepsilon} = \frac{\Pr(\tau < t) \cdot \Pr_b(|W_h + b - x| < \varepsilon)}{\Pr(|V_t - x| < \varepsilon)} = \frac{\Pr(|V_t - 2b + x| < \varepsilon)}{\Pr(|V_t - x| < \varepsilon)} \to \exp(-\frac{2b^2 - 2b(x+V_0) + 2xV_0}{\sigma^2 t})$$

as $\varepsilon \downarrow 0$

$$\Pr(V_A - C_0 > y | V_t = x) = \exp(-\frac{2(y+C_0)^2 - 2(y+C_0)(x+V_0) + 2xV_0}{\sigma^2 t})$$